\documentclass[11pt]{amsart}
\usepackage{amscd,amssymb}
\theoremstyle{plain}

\newtheorem{Thm}[equation]{Theorem}

\newtheorem{Prop}[equation]{Proposition}
\newtheorem{Lem}[equation]{Lemma}
\newtheorem{Def}[equation]{Definition}
\numberwithin{equation}{section}

\newcommand{\mc}[1]{{}}

\newcommand{\q}{\mathbb{Q}}

\newcommand{\e}{\epsilon}
\newcommand{\z}{\mathbb{Z}}
\renewcommand{\q}{\mathbb{Q}}
\newcommand{\n}{\mathbb{N}}

\newcommand{\br}{\mathbb{R}}

\newcommand{\ba}{\backslash}

\newcommand{\G}{\Gamma}

\newcommand{\ga}{\gamma}

\newcommand{\Cal}{\mathcal}
\renewcommand{\P}{\mathcal P}

\renewcommand{\gg}{\Gamma\ba G}
\newcommand{\la}{\langle}
\newcommand{\ra}{\rangle}
\newcommand{\pr}{\operatorname{pr}}

\newcommand{\lm}{\lambda_m}

\newcommand{\bi}{\begin{itemize}}

\newcommand{\ei}{\end{itemize}}

\newcommand{\SL}{\operatorname{SL}}

\newcommand{\GL}{\operatorname{GL}}

\newcommand{\Sp}{\operatorname{Sp}}

\newcommand{\vol}{\operatorname{vol}}

\begin{document}

\title[Representations of integers by an invariant polynomial]{
Representations of integers by an invariant polynomial and
 unipotent flows}
\author{Alex Eskin and Hee Oh}
\address{Mathematics Department, University of Chicago, Chicago,
IL 60637}

\thanks{The first author is partially supported by Packard
foundation.}
\thanks{The second author partially supported by NSF grants DMS 0070544
and DMS 0333397.}

\email{eskin@math.uchicago.edu}
\address{Mathematics Department, Princeton University, Princeton,
NJ 08544, Current address: Math 253-37, Caltech, Pasadena, CA 91125}
\email{heeoh@its.caltech.edu}
\maketitle
\section{Introduction}
 Let $f$ be an integral homogeneous polynomial of degree $d$
 in $n$ variables.
A basic problem in Diophantine analytic number theory is to
understand the behavior of the integral representations of
integers $m$ by $f$ as $m$ tends to infinity.

For each $m\in \n$, consider the level variety
 $V_m:=\{x\in \br ^n : f(x)=m\}$.
For instance, if $f(x_1, \cdots, x_n)=x_1^2+\cdots +x_n^2$ ($n\ge
3$), then $V_m$ is the sphere of radius $\sqrt m$ centered at the
origin and the set $V_m(\z)=V_m\cap \z^n$ consists of integral
vectors the sum of whose squares is equal to $m$. In this case,
the asymptotic of $\# V_m(\z)$ is well known. For $n\ge 5$ the
classical Hardy-Littlewood circle method applies and for $n= 4$
the Kloosterman sum method works. For $n=3$, Linnik gave a
conditional answer and later Iwaniec gave a complete answer (see
[Sa], [Iw]).

In the case when $V_m$ is non-compact, the number $\# V_m(\z)$
may be infinite. In this case one asks
if there exists an asymptotic density for $V_m(\z)$
as $m\to \infty$ . To be more precise, for a compact subset
$\Omega$ of $V_1$, set
 \begin{equation}\label{nmf}N_m(f,
\Omega):= \# V_m(\z)\cap m^{1/d}\Omega .\end{equation}

Or equivalently, $$N_m(f,
\Omega)
=\# V_m(\z)\cap \br ^+ \Omega $$ where $\br ^+\Omega$ is the
radial cone $\{x\in V: t x\in\Omega \text{ for some } t\in \br
^+\}$.
 The question we study in this paper is if there exists a sequence of
numbers $\omega_m$, independent of the compact subset $\Omega$,
such that
$$N_m(f, \Omega)\sim_{m\to \infty}  \omega_m \cdot \vol(\Omega).$$
This formulation is basically
 due to Linnik ([Li2], see also [Sa]).

 The only known general method is the Hardy-Littlewood
 circle method. However for this method to work,
 the number of variables needs to be much larger than the degree
 of the polynomial in general.

In this paper, we focus on the polynomials that are invariant
under an action of a semisimple real algebraic group. In such
cases, the level varieties admit actions of a semisimple algebraic
group and the dynamics of such groups plays a crucial role
in understanding this question.
 In particular, when $V_1$ is a
homogeneous space of a semisimple real algebraic group with the
stabilizer being generated by unipotent flows, we are able to use
a well developed theory of unipotent flows on a
homogeneous space of a Lie group.

We formulate our main results for a family of varieties $V_m$
which are more general than level sets of a given polynomial. Let
$V$ be a finite dimensional real vector space with a
$\z$-structure. Let $G$ be a semisimple real algebraic
  group defined over $\q$ and let
 $\rho:G \to \GL(V)$ be a $\q$-rational representation of $G$.
Fix a non-zero vector $v_0\in V(\z)$ such that $v_0G$ is Zariski
closed and
 a sequence $\{\lambda_m\in \br ^{+}
: m\in \n \}$ of strictly increasing numbers.
We set
\begin{equation*}\label{vm}
V_m:=\lm v_0 G \quad \text{ for each $m\in \n$} .\end{equation*}

Denote by $\Cal A_V$ be the collection of all arithmetic
subgroups $\G\subset G(\q) $ of $G$ such that
$V(\z)\G \subset V(\z)$. Since $\rho$ is defined over $\q$,
$\Cal A_V$ is non-empty. Consider
a sequence $\{\Cal O_m\subset V_m(\z): m\in \n\}$
of $\G$-invariant subsets of $V_m(\z)$
for some $\G\in \Cal A_V$.
For a compact subset
$\Omega$ of $V_1$, we define
$$N(\Cal O_m, \Omega):= \# \operatorname{pr}(\Cal O_m )\cap \Omega$$
where $\operatorname{pr}:V_m \to V_1$ denotes the radial projection given by
$\pr(x)= \lambda_1 \lm ^{-1}x$.

 Let $H$ denote the stabilizer of
$v_0$ in $G$. The notation $H^0$ and $G^0$ denote the identity
components of $H$ and $G$ respectively. In the rest of the
introduction, we
 assume that $H$ is semisimple without compact
factors.

\begin{Thm}\label{maximal} Suppose that $H^0$ is a
 maximal connected closed subgroup of $G^0$.
Let $\{\Cal O_m\subset V_m(\z)\}$ be a sequence of $\G$-invariant
subsets for some $\G\in \Cal A_V$, for example, $\Cal O_m=V_m(\z)$.
 Suppose that for each $m_0\in
\n$,
\begin{equation}\label{prim} \# \{m\in \n: \pr (\Cal O_{m_0})=
\pr (\Cal O_{m})\}< \infty .\end{equation} Fix a $G$-invariant
Borel measure $\mu$ on $V_1$.
 Then for any compact subset
$\Omega$ of $V_1$ with boundary of measure $0$,
 \begin{equation}\label{p1} N(\Cal O_m,  \Omega)\sim_{m \to \infty}
  \omega (\Cal O _m) \cdot \mu{(\Omega )} \end{equation}
where $\omega (\Cal O_m)$ is given below (\ref{omega}).
\end{Thm}

An immediate consequence of the above theorem
is
that the sequence $\{\pr(\Cal O_m)\}$ of projections is
equidistributed on $V_1$ as $m\to \infty$,
since it follows from $(\ref{p1})$ that
for any two compact subsets $\Omega_1$ and $\Omega_2$
of $V_1$ with boundary measure $0$,
$$\frac{N_m(\Cal O_m, \Omega_1)}{N_m(\Cal O_m, \Omega_2)} \sim_{m \to \infty}
  \frac {\mu{(\Omega_1 )}}{\mu(\Omega_2)}. $$

The condition (\ref{prim}) says that there is no infinite sequence
of varieties $V_{m_i}$ where all the integer points in $\Cal
O_{m_i}$ are simply the radial dilations of a fixed $\Cal O_{m_1}$.
Clearly this is a necessary condition for the conclusion of the
above theorem. We remark that the condition (\ref{prim}) is also
equivalent to saying $\lim _{m\to \infty}\omega(\Cal O_m)=\infty$
(see Lemma \ref {q}).

In order to give the formula of the number $\omega(\Cal
O_m)$ in Theorem \ref{maximal}, let $\mu_G$ and $\mu_{H}$ denote
Haar measures on $G$ and $H$ respectively so that the triple
$(\mu, \mu_G, \mu_H)$ are topologically compatible in the sense of
Weil [We], namely, for any continuous function $f$ on $G$ with
compact support,
\begin{equation}\label{compatible} \int_G f\, d\mu_G(g)=
\int_{H\ba G} d\mu(\bar g) \int_H f(hg) d\mu_H(h)
\end{equation}
 where $\mu$, by slight abuse of notation, denotes the
$G$-invariant measure on $H\ba G$ pulled back from $\mu$ on $V_1$ via
the identification of $H\ba G$ with $V_1$ by $Hg\mapsto v_0g$.

The number $\omega (\Cal O_m)$ is defined as follows: for any $\G\in
\Cal A_V$ preserving each $\Cal O_m$,
\begin{equation}\label{omega}
\omega (\Cal O_m):=  \frac{\sum _{\xi\G\subset \Cal O_m} \mu_H
\left( ( H\cap {g_\xi}^{-1} \G g_{\xi}) \ba H\right)}{\mu_G(\G \ba
G)},\end{equation} where the sum is taken over the set of disjoint
$\G$-orbits in $\Cal O_m$ and $g_\xi\in G$ is an element such that
$v_0=\operatorname{pr}(\xi)g_\xi$.
 Under our assumption, $\omega (\Cal O_m) <\infty $, and moreover
  $\omega(\Cal O_m)$ does not depend on the choice of $\G$ in
  $\Cal A_V$ preserving $\Cal O_m$, justifying our notation.
 This is shown in [Oh]
  for $\Cal O_m=V_m(\z)$ but the same proof works for
  any $\Cal O_m$. For the case of $\Cal O_m=V_m(\z)$,
  we simply write
  \begin{equation}\label{wm}\omega_m=\omega(V_m(\z)).\end{equation}

 Note also that the product $\omega (\Cal O_m)
\cdot \mu(\Omega)$ is independent of the choice of compatible
measures $(\mu,
\mu_H, \mu_G)$.
For this reason, we sometimes write $\omega(\Cal O_m)\cdot
\operatorname{vol}(\Omega)$ in what follows.

 We remark that Theorem \ref{maximal} applies to the cases where $V_1$
is an affine symmetric space, i.e., $H$ is the set of fixed points
of an involution of $G$, with additional assumptions that $G$ is
$\q$-simple and that $H$ is semisimple without compact factors. It
is so since $H^0$ is then a maximal connected closed subgroup of
$G^0$ (cf. [Bo, Lemma 8.0]).

We present some examples which follow
from Theorem \ref{maximal}
 for the case $\Cal O_m=V_m(\z)$. See (\ref{nmf})
for the notation $N_m(f, \Omega)$. The following three theorems are
proven in section \ref{example} where we realize each $f$ as an
invariant polynomial of a certain representation of a semisimple
algebraic group. Once we do that, the number $\omega_m$ is defined
as in (\ref{omega}) and (\ref{wm}).

 \begin{Thm}\label{determinant} Let $n \ge 2$
 and $\operatorname{Det}$ denote the determinant polynomial on
the vector space $\operatorname{M}_n(\br)$ of $n\times n$ matrices.
For any compact subset $\Omega$ of $V_1=\{X\in
\operatorname{M}_n(\br): \operatorname{Det}(X)=1\}$ of boundary of
measure $0$, we have
\begin{equation*}
N_m(\operatorname{Det}, \Omega) \sim_{m\to \infty} \omega_m \cdot
\operatorname{vol}{(\Omega)} .
\end{equation*}
\end{Thm}

Theorem \ref{determinant} was proven first by Linnik
[Li1] (for $n=2$) and by Linnik and Skubenko (for $n \ge 3$) [LS] in the early
sixties using methods in analytic number theory. A different proof
was also given using methods based on Hecke operators (see [Sa],
[COU] and [GO]).

Let $V$ be the subspace of $\operatorname{M}_{2n}(\br)$ consisting
of the skew-symmetric matrices, i.e., $V=\{X\in
\operatorname{M_{2n}}(\br) : X^t=-X\}$ ($n\ge 2$) and consider the
Pfaffian as a polynomial on $V$ defined by $\operatorname{Pff}(X)
^2=\operatorname{Det}(X)$ and
$$\operatorname{Pff}\begin{pmatrix}
0 &
I_n \\
-I_n& 0\end{pmatrix} =1$$ where $I_n$ denotes the identity matrix of order
$n$.

 \begin{Thm}\label{pfaffian} Let $n \ge 2$.
For any  compact subset $\Omega$ of $V_1=\{X\in V:
\operatorname{Pff}(X)=1\}$ of boundary of measure $0$, we have
\begin{equation*}
N_m(\operatorname{Pff}, \Omega) \sim_{m\to \infty} \omega_m \cdot
\operatorname{vol}{(\Omega)}.
\end{equation*}
\end{Thm}

 An integer $m$ is called a
fundamental discriminant if and only if $m$ is either a square
free integer congruent to $1$ mod $4$, or $4$ times of a square
free integer which is congruent to $2$ or $3$ mod $4$.
\begin{Thm}\label{quadratic}
Let $Q$ be an integral quadratic form of signature $(r, s)$ where
$r+s\geq 4$, $r\ge 2$ and $ s\ge 1$.
 For any compact subset $\Omega$ of
$ V_1=\{X\in \br^{r+s}: Q(X)=1\}$ of boundary of measure $0$,
 $$N_m(Q, \Omega)
\sim
 \omega_m \cdot \operatorname{vol}{(\Omega)}$$
as $m \to \infty$ along the fundamental discriminants.
 \end{Thm}
In fact, the Hardy-Littlewood circle method together with the
Kloosterman sum method (needed for $r+s=4$) (cf. [Va], [Es]) also
gives an asymptotic density in Theorem \ref{quadratic}, in the
form of product of local densities. By comparing the two different
forms of the asymptotic for $N_m(Q, \Omega)$, one obtains a new
proof of Siegel mass formula for quadratic forms (see [Oh]) in the
same spirit of the work of Eskin, Rudnick and Sarnak [ERS].
The above theorem is still true for the case of $r+s=3$ and $rs>0$
though our method does not apply. This follows from a theorem of Duke [Du].

More generally, a natural question is whether the asymptotic $\omega_m\cdot
\operatorname{vol}(\Omega)$ in Theorem \ref{maximal} coincides with the
heuristics predicted by the Hardy-Littlewood circle method. It is
shown in [Oh] that this is true in many cases but not always.

 To state our main result in a more general setting without
  the maximality assumption on $H$, we define the following:
the notation $Z(H^0)$ denotes the centralizer of $H^0$ in $G$.

\begin{Def}\label{foc} A sequence $\{\Cal O_m\subset V_m(\z)\}$
 is called
{\bf focused} if, for every $\G\in \Cal A_V$ preserving $\{\Cal
O_m\}$, there exist a proper connected closed subgroup $L$ of $ G^0$ in
which $L\cap \G $ is a Zariski dense lattice, $g\in G$ with
$g^{-1}H^0 g\subset L$ and a compact subset $\Cal C\subset v_0
Z(H^0)gL$ such that
 for every compact subset $\Omega \subset V_1$, there exists
$\gamma_\Omega\in \G$ with
\begin{equation}\label{prp} \limsup_{m \to
\infty}\frac{\sum_{\xi\G\subset \Cal O_m}
 \{\omega (\xi \G) :\operatorname{pr}
(\xi \G)\cap \Omega\subset \Cal C (L \cap \G)
\gamma_\Omega\}}{\omega (\Cal O_m)} >0 . \end{equation}
\end{Def}
In the above and also in the rest of the paper,
 the notation $\sum_{\xi\G\subset \Cal O_m}$
means that the sum is taken over the set of disjoint $\G$-orbits
contained in $\Cal O_m$.

\begin{Thm}\label{Main}\label{main} Suppose that
 $\{\Cal O_m: m\in \n\}$
 is not focused. Then
for any compact subset $\Omega$ of $ V_1$ with boundary of measure
$0$,
 $$N(\Cal O_m , \Omega)\sim_{m \to \infty}
  \omega (\Cal O_m) \cdot \operatorname{vol}(\Omega) .$$
\end{Thm}

Note that in the case when $\Cal O_m$ consists of one $\G$-orbit,
the focusing of $\{\Cal O_m\}$ implies that for every given compact
subset $\Omega$ of $V_1$, there exists an infinite sequence $m_i$
such that $\pr(\Cal O_{m_i})\cap \Omega $ lies completely inside a
proper subvariety of $V_1$.
 In general, if $\{\Cal O_m\}$ is focused, one
expects in view of Theorem \ref{main} that a positive proportion of
points in $\pr(\Cal O_{m_i})\cap \Omega$ lies inside a proper
subvariety of $V_1$ of the form $v_0Z(H)g L\gamma$ for some
$\gamma\in \G$ and for some proper subgroup $L$ of $G^0$.

Assuming $H$ connected for simplicity, we explain some of schemes in
the proofs. For each $\xi\G\subset \Cal O_m$, consider the
$H$-invariant measure $\nu_{\xi\G}$
 on $\G\ba G$ supported on the closed orbit $\G\ba \G g_\xi H$
with the total measure $\omega (\xi\G)$. Let $\sigma_m$ denote the
averaging measure $\frac{\sum_{\xi\G\subset\Cal
O_m}\nu_{\xi\G}}{\omega(\Cal O_m)}$ on $\G\ba G$. We first observe
that the equidistribution of $\pr(\Cal O_m)$ (Theorem \ref{main})
follows if $\sigma_m$ converges to the $G$-invariant probability
measure on $\G \ba G$ as $m \to \infty$ (Proposition \ref{Av}). We
then show that under the non-focusing assumption of $\{\Cal O_m\}$,
the $H$-orbits $\G\ba \G g_\xi H$ which stay outside a given compact
subset $C$ as well as the $H$-orbits $\G\ba \G g_\xi H$ which stay
inside $\G\ba \G Lg_\xi$ for some proper subgroup $L$ of $G$ are ignorable,
in considering the weak limits of $\sigma_m$. Then by applying main
ergodic results of Dani and Margulis in [DM1-2], we show that each $H$-orbit
$\G\ba \G g_\xi H$ is getting longer and longer and moreover
uniformly distributed on $\G \ba G$ as $m \to \infty$; hence the
average measure $\sigma_m$ tends to the Haar measure on $\G \ba G$
as $m \to \infty$.

In particular in the case where the
 subgroup
$H$ is a maximal closed subgroup
of $G^0$, the equidistribution of $\{\pr (\Cal O_m)\}$ is
a consequence of the phenomenon that for any sequence
of non-repeated individual $\G$-orbits $\{\pr (\xi) \G \subset
\pr(\Cal O_m)\}$,
 the corresponding sequence $\{\G\ba\G g_\xi H\}$ of $H$-orbits is
 uniformly distributed on $\G \ba G$ and hence
 $\pr (\xi) \G$ is equidistributed on $V_1$ as $m \to \infty$.

 \noindent{\bf Remark }   Gan and Oh [GO] showed
 that for any invariant polynomial $f$ in the above setting
but with a more general $H$ (not necessarily semisimple), if
$V_{d_0}(\z)$ contains at least one integer point for some $d_0\in
\n$, then there exist explicitly computable constants $c$ and $r$
depending only on $G$, $\rho$ and $\text{deg }(f)$ such that
$\operatorname{pr}(V_{cd_0m^r}(\z))$ becomes dense on $V_1$ in a
strong sense as $m \to \infty$. In particular it follows that the
limit supreme in
 (\ref{prp}) is strictly less than
$1$.

\noindent{\bf Acknowledgment}
The second named author would like to thank Peter Sarnak for
suggesting this problem
as well as for helpful discussions.

\section{\bf Measure theoretic formulation of a counting
problem}\label{count} Let $G$ be a real semisimple
 algebraic group defined over $\q$. This means that
there exists a connected semisimple algebraic group $\underline G$
defined over $\q$ such that $G$ is a closed subgroup of $\underline
G(\br)$  containing the identity component $\underline G(\br)
^\circ$. Let $\rho:G \to \GL(V) $ be a $\q$-rational representation
for a finite dimensional real vector space $V$ defined over $\z$.
Fix a non-zero vector $v_0\in V$ such that $v_0G$ is Zariski closed
and a sequence $\{\lm : m\in \n\}$ of strictly increasing positive
numbers. For each $m\in \n$, we set
\begin{equation*} V_m=\lm v_0 G . \end{equation*}
Let $H$ denote the stabilizer of $v_0$ in $G$.
Since $\rho$ is rational, $H$ is a real algebraic subgroup of $G$.
In particular, the identity component $H^0$ is a finite index normal
subgroup of $H$.
 Assume that $H^0$ has no non-trivial $\br$-character.
 Let $\mu_{G}$ denote
 a Haar measure on $G$.
Since both $G$ and $H$ are unimodular, we may choose a
$G$-invariant Borel measure $\mu$ on $V_1$ and an $H$-invariant
measure $\mu_H$
 on $H$ so that the triple $(\mu, \mu_G,\mu_H)$ is compatible
 in the sense of (\ref{compatible}).

Let $\Cal A_V$ be defined as in the introduction, and let $\G\in \Cal A_V$.
Since $G$ acts transitively on
$V_1$, for any $\xi\in V_m$,
 there exists $g_\xi\in G$  such that $v_0=\pr (\xi) g_\xi$.
The choice of $g_\xi$ is unique only up to modulo $H$.
If $H_\xi$ denotes the stabilizer of $\xi$ in $G$
for $\xi\in V_m(\z)$, then $H_\xi$ is a $\q$-subgroup of $G$
and $H_\xi=g_\xi Hg_\xi^{-1}$. Therefore $H_\xi^0$
has no non-trivial $\q$-character and hence by a theorem of
Borel and Harish-Chandra, $H_\xi\cap \G$ is a lattice in $H_\xi$.
Hence
$$\omega(\xi\G):=
\frac{\mu_H(H\cap g_\xi^{-1}\G g_\xi\ba H)}{\mu_G(\G\ba G)} <\infty .$$
 Observe that the definition of
$\omega(\xi\G)$ depends only on the $\G$-orbit not on its
representative.

Let $\{\Cal O_m\subset V_m(\z)\}$ be
a sequence of non-empty $\G$-invariant subsets of $V_m(\z)$
 for some $\G\in \Cal A_V$.
  Since each $V_m$ is Zariski closed, by a theorem of
Borel and Harish-Chandra, the number of $\G$-orbits in
 $V_m(\z)$ and hence in $\Cal O_m$ is finite.
Hence
$$\omega(\Cal O_m):=
\sum_{\xi\G\subset \Cal O_m} \omega(\xi\G) <\infty .$$

\begin{Lem}[Oh]
The number $\omega(\Cal O_m)$ is
independent of the choice of $\G\in \Cal A_V$
preserving $\Cal O_m$.
\end{Lem}

The space $\P (\G \ba G)$ of the probability measures on $\G \ba G$
is equipped with the weak$^*$-topology.
Now fix any $\G\in\Cal A_V$ which preserves each $\Cal O_m$.
For each  $\G$-orbit $\xi\G \subset\Cal O_m$,
let $\nu_{\xi\G }$ denote the unique $H$-invariant measure
on $\G\ba G$ supported on the closed orbit $\G\ba \G g_\xi H$
and with the total measure given by
$\omega(\xi\G)$.
Hence $$\frac{1}{\omega(\Cal O_m)}
\sum_{\xi\G\subset \Cal O_m}\nu_{\xi\G}\in \P(\G\ba G) .$$

Here is a main proposition suggested by Sarnak
which translates the counting problem
 to the question of whether the weak-limits of the above
 measures are $G$-invariant.

\begin{Prop}\label{Av}
If $$\lim_{m\to
\infty}\frac{1}{ \omega (\Cal O_m)} \sum_{\xi\G\subset \Cal O_m}\nu_{\xi\G}
 =\frac{1}{\mu_G(\G\ba G)} \mu_{ G} \quad \text{ in $\P (\G \ba G)$}$$
 then
for any compact subset $\Omega$ of $ V_1$ with boundary of measure $0$,
$$ N(\Cal O_m, \Omega) \sim_{m \to \infty}
{\omega (\Cal O_m)} \cdot \mu{(\Omega )}.$$
\end{Prop}
\begin{proof} Without loss of generality, we assume $\mu_G(\G\ba G)=1$.
Let $\phi$ be any continuous
function with compact support
on $H\ba G =V_1$.
Define a function $F^m_{\phi}$ as follows: for each $g\in G$
$$F ^m_{\phi}(g):=\frac{1}{\omega (\Cal O_m)} \sum_{\xi\G\subset \Cal O_m}
\sum_{\gamma \in (H_{\xi}\cap \G) \ba \G} \phi  (\pr(\xi) \ga g).$$
Since $F^m_{\phi}$ is left $\G$-invariant,
it may be considered as a function on $\G\ba G$.
Let $\psi $ be a continuous
function on $\G \ba G$ with compact support.
Note that

 \begin{align}
 {\omega (\Cal O_m)}\cdot \la F^m_{\phi} , \psi \ra
&=
\sum_{\xi\G\subset \Cal O_m} \int_{\gg}
\left(
\sum_{\gamma \in (H_{\xi}\cap \G) \ba \G} \phi  (\pr(\xi)\ga g)\psi (g)\right)
\, d\mu_G(g)\notag \\
&=
\sum_{\xi\G\subset \Cal O_m} \int _{g\in (H_{\xi}\cap \G)\ba
G} \phi (\pr(\xi) g) \psi (g)\, d\mu_{G}(g) \notag \\
&= \sum_{\xi\G\subset \Cal O_m}
\int _{t\in (H\cap g_\xi^{-1}\G g_\xi) \ba
G} \phi (\pr(\xi) g_\xi t ) \, \psi (g_\xi t)\, d\mu_{G}(t) \notag \\
&=  \sum_{\xi\G\subset \Cal O_m}
\int_{g\in H\ba G}\phi (v_0 g ) \left(  \int _{h\in (H\cap g_\xi^{-1}\G g_\xi )  \ba H}
  \psi (g_\xi hg ) d\mu_H(h) \right) \, d\mu (g) \notag\\
&=  \sum_{\xi\G\subset \Cal O_m}
\int_{g\in H\ba G}\phi (v_0g ) \left( \int _{s\in \G \ba G}
  \psi (sg )\, d\nu_\xi (s) \right)\, d\mu (g)\notag
\end{align}

Consider a function $\psi_g$ on $\G \ba G$
 defined by $\psi_g(s):=\psi(sg)$. Then
$$ \int_{\G \ba G} \psi_g \, d\mu_G= \int_{\G \ba G} \psi \, d\mu_G .$$
Hence by the assumption,
$$
\lim_{m\to \infty} \frac{1}{\omega(\Cal O_m)} \sum_{\xi\G\subset \Cal O_m}
 \int _{s\in \G \ba G}
  \psi (sg )\, d\nu_\xi (s)
= \int_{\G \ba G} \psi \, d\mu_G .$$

Now by the
 Lebesgue dominated convergence theorem,
$$
\lim_{m\to \infty} \la  F_{\phi} ^m, \psi \ra
=\int_{H\ba G}\phi (v_0g)  d\mu(g)
\cdot  \int_{\G \ba G} \psi \, d\mu_G.$$
If follows that
\begin{equation}\label{par}
\lim_{m \to \infty}
\la  F_{\chi_{\Omega}} ^m, \psi \ra = \mu{(\Omega)}
\cdot  \int_{\G \ba G} \psi \, d\mu_G \end{equation}
where $\chi_{\Omega}$ denotes the characteristic function of $\Omega$.

Fix $\e>0$.
Let $U_\e$ be a symmetric neighborhood of $e$ in $G$ such that
$$\mu(\Omega_\e^+-\Omega_\e^-)\le \e$$
where $\Omega_\e^+=\cup_{u\in U_\e} \Omega u$ and $\Omega_\e^-=
\cap_{u\in U_\e}\Omega u$.
Then for all $g\in U_{\e}$,
\begin{equation}\label{Fi}
F_{\chi_{\Omega _{\e -}}} ^m (g) \leq
F_{\chi_{\Omega}} ^m (e) \leq F_{\chi_{\Omega _{\e +}}} ^m (g).\end{equation}

Let $\psi_{\epsilon}$ be
a non-negative continuous function on $\G\ba G$ with support in
$U_{\e}$  and $\int_{\G\ba G}\psi_\e \, d\mu_G =1$.
Integrating (\ref{Fi}) against
$\psi_\e$ now gives
$$\la F_{\chi_{\Omega _{\e -}}} ^m , \psi_{\e}\ra  \leq
F_{\chi_{\Omega}} ^m (e) \leq \la F_{\chi_{\Omega _{\e +}}} ^m,
 \psi _{\e} \ra .$$
Since both sides tend to $\mu (\Omega _{\e \pm})$
respectively as $m \to \infty$ by (\ref{par}) and $\e >0$ is
arbitrary, we have $$F_{\chi_{\Omega}} ^m (e)\to \mu(\Omega)
\quad \text{ as $m \to \infty$}.$$  Since
$$F^m_{\chi_{\Omega}}(e)
 =
\frac{N(\Cal O_m,  \Omega) }{ \omega (\Cal O_m) },$$
this proves the claim.
\end{proof}

\section{\bf Asymptotic behavior of unipotent flows}
We recall the following
 fundamental result of Dani and Margulis.

\begin{Thm} [DM2, Theorem 6.1]\label{di} Let $G$ be a connected Lie group and
$\G$ a lattice in $G$. Given a compact subset $C\subset \G\ba G$
and an $\e >0$, there exists a compact set $K\subset \G\ba G$ such
that the following holds: for any $x\in C$, any unipotent
one-parameter
 subgroup $\{u(t)\}$ of $G$, and any $T>0$,
$$|\{t\in [0, T]:xu(t)\in K\}| > (1-\e)T$$
where $|\cdot |$ denotes the Lebesgue measure on $\br$.
\end{Thm}

Let $G$ be a connected semisimple real algebraic group defined over $\q$, and
 $H$ a connected semisimple real algebraic subgroup of $G$.
 Let $\G\subset G(\q)$ be an arithmetic subgroup of $G$. The results
in this section have meanings only when $\G\ba G$ is non-compact, which
we assume. Consider the
one point compactification $\G\ba G \cup \{\infty\}$ of $\G\ba G$. The space $ \P
(\G\ba G\cup\{\infty\})$ of the probability measures on
$\G\ba G\cup\{\infty\}$ equipped with the weak$^*$-topology is weak$^*$
compact.

 Let $\{g_m \in G\}$ be a sequence such that $g_m H g_m^{-1}$ is
a $\q$-subgroup of $G$ for each $m$.
 By a theorem of Borel and Harish-Chandra [BH],
it follows that $g_m^{-1}\G g_m \cap H$ is a lattice in $H$.
Hence each $\G\ba \G g_mH$ is closed in $\G\ba G$ (cf. [Rag])
and there exists the unique $H$-invariant
probability measure $\mu_m$ in $\G\ba G$ supported on $\G\ba \G g_m H$.

\begin {Prop}\label{supp0}
Assume either that $H$ has no compact factors
or that $g_mHg_m^{-1}$ is $\q$-simple for each $m$.
Then the following are equivalent:
\begin{itemize}
\item[(1)] There exists a compact subset $C$ of $ \G\ba G$
such that
$$\G\ba \G g_m H \cap C \ne \emptyset \quad\text{for all
sufficiently large $m \in \n$.}$$
\item[(2)] Every weak limit of $\{\mu_m\}$ in $\P(\G\ba G\cup\{\infty\})$
 is supported on $\G\ba G$.
\end{itemize}
 \end{Prop}
\begin{proof} Assume that (1) is true.
 Without loss of generality, we may assume that
$\G\ba \G g_m \in C$ for all $m \in \n$. Let $H_N$ denote
the unique maximal connected normal closed subgroup of $H$ without compact factors.
Let $U=\{u(t)\}$ be a
unipotent one parameter subgroup in $H_N$ not contained in any
proper normal subgroup of $H_N$.
Such a subgroup exists (see for example, [MS, Lemma 2.3]).
Under our assumption, either $H=H_N$ or $g_m^{-1}\G g_m\cap H$
is an irreducible lattice in $H$. Hence it follows from Moore's ergodicity
theorem
 (cf. Theorem 2.1 in [BM]) that
$U$ acts ergodically with respect to each $\mu_m$.
Moreover by
 the Birkhoff ergodic theorem (cf. [BM]), for almost all $h\in H$,
$\G\ba \G g_m hu(t)$ is uniformly distributed on $\G\ba G$ with respect
to $\mu_m$, i.e., for any $f\in C_c(\G\ba G)$,
$$\lim_{T\to \infty} \frac{1}{T}\int_0^Tf(\G\ba \G g_m h u(t) )\, dt
=\int f \, d\mu_m .$$
Therefore we may assume that
for each $m \in \n$, there exists $h_m \in H$ such that $\G\ba \G g_m h_m \in C$
and $\G g_m h_mu(t)$
 is uniformly distributed on $\G\ba G$ with respect to
$\mu_m$. For any
given $\e >0$, let $K$ be a compact subset of $\G\ba G$ as in Theorem
\ref{di} with respect to $C$. Then for each $m\in \n$, $$\mu_m (K) \ge 1-\e .$$
Therefore $\mu (K)\ge 1-\e$ for any weak limit $\mu$ of $\{\mu_m\}$. Since
$\e>0$ is arbitrary, we have $\mu (\G\ba G)=1$, proving that (2) holds.

Now suppose that (1) fails.
First write $\G\ba G$ as $\cup_{i=1}^{\infty} C_i$
where $C_i$ are compact subsets such that
$C_i\subset C_{i+1}$ for all $i$. Then for each $i$, there exists
$m_i$ such that $$\G\ba \G g_{m_i}H\cap C_i=\emptyset .$$
Since $C_i$ is increasing,
we have
\begin{equation}\label{ci}
\G\ba \G g_{m_j}H\cap C_i=\emptyset\quad  \text{ for all $j\ge
i$}.\end{equation} This implies that any weak limit of
$\{\mu_{m_i}\}$ cannot be supported on $\G\ba G$, for if so,
  then for some $i_0$, $\mu_{m_i}(C_{i_0})>1/2$ for infinitely many $i$.
This is contradiction to (\ref{ci}).
 Hence (2) implies (1).
 \end{proof}

Let $\{\frak O_m\subset H\ba G\}$ be a sequence of
a finitely many union of $\G$-orbits.
For each $H\ba Hg\in \frak O_m$,
we assume that $gHg^{-1}$ is a $\q$-subgroup of $G$.
For each $\G$-orbit $\eta \G\subset\frak O_m$, set
$$\omega(\eta \G):=\frac{\mu_H((H\cap g_\eta ^{-1}\G g_\eta)\ba H)}{\mu_G(\G\ba G)} \quad\text{ and} \quad \omega(\frak O_m):=
\sum_{\eta\G\subset\frak O_m}\omega (\eta\G) .$$
Here $g_\eta\in G$ is such that $\eta=H\ba Hg_\eta$.

Let $\nu_{\eta\G}$ denote the $H$-invariant measure on $\G\ba G$ supported
on $\G\ba \G g_{\eta}H$ with the total measure given by
$\omega(\eta\G)$.
 Define an $H$-invariant probability measure $\sigma_m$ on $\G\ba G$:
$$\sigma_m=\frac 1{\omega(\frak O_m)} \sum_{\eta\G\subset \frak O_m}
\nu_{\eta \G}.$$

The notation $H_N$ denotes the unique maximal connected normal closed
subgroup
of $H$ without compact factors.
\begin{Prop}\label{Non-div}
Assume either that $H$ has no compact factors or
that $g_\eta H g_\eta^{-1}$ is $\q$-simple for any $\eta\in \frak O_m$.
Suppose that
$g_\eta H_N g_\eta^{-1}$
 is not contained in any proper $\q$-parabolic subgroup of $G$
for any $\eta\in \cup_m \frak O_m$. Then any weak
limit of $\{\sigma_{m} :m\in \n\}$ in $\P (\G\ba G\cup
\{\infty\})$ is supported in $\G\ba G$.
\end{Prop} \begin{proof}
Without loss of generality we may assume that
$\{\sigma_{m}\}$ converges in $\P (\G\ba G\cup\{\infty\})$. It
suffices to show that,
 for any $\e >0$,
there exists
 a compact subset $K\subset \G\ba G$ such that $$\sigma_{m}(K)>1-\e
\quad\text{for all sufficiently large $m$.}$$

Assume not; then for any compact subset $K\subset \G\ba G$, (after
going to a subsequence) there exists $\eta_m\in\frak O_m$ such that
$$\nu_{\eta_m\G} (K)<(1-\e) \omega(\eta_m\G) \quad\text{ for each $m\in \n$.}$$
  Let
 $U=\{u(t)\}$ be a unipotent one-parameter
subgroup of $H_N$ as in the proof of previous proposition. Let $R$
be the set of $h\in H$ such that $\G\ba \G g_{\eta_m} h u(t)$ is
uniformly distributed in $\G\ba \G g_{\eta_m}H$ with respect to the
probability measure $\frac 1 {\omega(\eta_m \G)} \nu_{\eta_m\G}$.
Then $R$ has the full measure in $H$ (see the proof of the
previous proposition).
 Fix any $h\in R$. Then for each $m\in \n$, there exists
$T_m>0$ (depending on $h$) such that
$$\frac{1}{T}|\{t\in [0, T]: \G\ba \G g_{\eta_m} hu(t)\in K\}| <1- \e/2 $$
for all $T>T_m$ where $| \cdot |$ denotes the
Lebesgue measure on $\mathbb R$. Applying a theorem of Dani and Margulis
[DM1, Theorem 2] (see also [EMS1]), we obtain that
for any
given $\alpha _m>0$ with $\lim_{m\to \infty}\alpha_m =0$,
 after passing to a subsequence,
there exist a proper parabolic $\q$-subgroup
$P$ of $G$, a non-zero vector $q\in \wedge ^k \operatorname{Lie} (W)(\q)$
 ($W$ being the unipotent radical of $P$ and $k=\operatorname{dim}(W)$)
and a sequence $\{\ga_m(h)\in \G\}$ such that for
all $m\in \n$ and $t>0$,
 $$\| q. \ga _m(h) g_{\eta_m} hu(t)\| <\alpha _m$$ where the
action is through the $k$-th exterior of the adjoint representation of $G$
 on $\wedge ^k \operatorname{Lie} (G)$.
Since $u(t)$ acts as a unipotent
one-parameter subgroup on $\wedge ^n \operatorname{Lie} (G)$
and any orbit of a unipotent one-parameter subgroup is unbounded except
for a fixed point,
it follows that for all $0\le t<\infty $ and
for all $m\in \n$,
 $$q.\ga_m^h g_{\eta_m} hu(t)=q.\ga_m^h g_{\eta_m} h .$$
Hence
 $$U \subset {(\ga_m^h g_{\eta_m} h)}^{-1}P (\ga_m^h
  g_{\eta_m} h),
$$ since the latter group contains
 the stabilizer of the vector $ q.\ga_m^h g_{\eta_m} h$.
Hence we have shown that for almost all $h\in H$ and for any $m\in
\n$,
$$hUh^{-1}\subset {(\ga_m^h g_{\eta_m})}^{-1} P(\ga_m^h g_{\eta_m})$$
 for some $\gamma_m^h\in \G$.
Since $\G$ is countable, it follows that for each $m\in \n$, there
exist an element $\ga_m\in \G$ and a subset $S_m\subset H$ of
positive measure such that $$hUh^{-1}\subset {(\ga_m
g_{\eta_m})}^{-1} P(\ga_m g_{\eta_m})\quad\text{for all $h\in S_m$} .$$ Since
the set
$$\{h\in H: hUh^{-1}\subset {(\ga_m g_{\eta_m})}^{-1} P(\ga_m
g_{\eta_m})\}$$ is a real analytic submanifold of $H$ with a positive measure,
 it is indeed equal to $H$. Hence
$$hUh^{-1}\subset {(\ga_m g_{\eta_m})}^{-1} P(\ga_m g_{\eta_m})\quad
\text{ for all $h\in H$.}$$
Since $U$ is not contained in any proper normal subgroup of $H_N$,
it follows that
$$g_{\eta_m}H_Ng_{\eta_m}^{-1} \subset {\ga_m}^{-1} P\ga_m .$$
This contradicts the assumption since ${\ga_m}^{-1} P\ga_m$ is a proper
parabolic
$\q$-subgroup of $G$.
\end{proof}

\section{Projections of $\Cal O_m$ and stabilizer subgroups}
We recall the following theorem of Dani and Margulis: let $G$ be
any connected Lie group and $\G$ a discrete subgroup of $G$. We
fix a left invariant Riemannian metric on $G$. Let
$M$ be any closed subgroup of $G$ such that $M\cap \G$ is a
lattice in $M$. Then $\G\ba \G M$ is a closed Riemannian
submanifold
of $\G\ba G$ and hence it has a
right $M$-invariant
 Riemannian volume form, denoted by $\mathcal V$,
induced by the Riemannian metric.

\begin{Thm}[DM2, Theorem 5.1]\label{fv}
For any $c >0$, let $\Cal W_c$ be the collection of all closed
connected
 subgroups of $G$ such that $\G\ba \G M$ is closed in $\gg$ and
 $\mathcal V((M\cap \G) \ba M)\le c$.
Then there are only finitely many subgroups of the form $M\cap \G$
 with $M\in \Cal W_c$.
 \end{Thm}

We also need the following simple consequence of a theorem of
Kazhdan and Margulis ([KM], [Ra, Theorem 11.8]):
\begin{Lem}\label{km}
Let $G$ be a connected linear semisimple Lie group without compact
factors.
 There exists a constant
$c>0$ such that for any discrete subgroup $\G$ of $G$, the
co-volume of $\G$ in $G$ (with respect to a fixed Haar measure
on $G$) is at least $c$.
\end{Lem}

We keep the same notation from section \ref{count} for
 $G$, $\rho:G \to \GL(V)$, $\G$, $H$, $H_\xi$, $v_0$, etc.
 Let $\Cal O_m$ be a
$\G$-invariant subset of $V_m(\z)$ for each $m$.
We assume that $H$ is semisimple without compact factors.
Denote by $N(H)$
the normalizer of $H$ in $G$.
\begin{Lem}\label{q} Assume that $[N(H):H] < \infty$ and
that
for each $\xi\in \Cal O_m$, $H_\xi^0$ is not contained in any proper
$\q$-parabolic subgroup of $G$. Then the following are equivalent:
\begin{itemize}
\item[(1)] Suppose that for each $m_0\in
\n$,
\begin{equation*} \# \{m\in \n: \pr (\Cal O_{m_0})=
\pr (\Cal O_{m})\}< \infty .\end{equation*}
\item[(2)] $\lim_{m\to \infty}\omega(\Cal O_m)= \infty$.
\end{itemize}
  \end{Lem}

\begin{proof}
It is easy to see that $\pr(\Cal O_m)=\pr (\Cal O_k)$ implies $\omega(\Cal
O_m)=\omega(\Cal O_k)$. Hence (2) implies (1).
Assume now that (2) fails. Then by passing to a subsequence we may
assume that $\omega(\Cal O_m)$ is uniformly bounded.
By Lemma \ref{km} which we may apply since $H^0$
has finite index in $H$, there exists some $c>0$ such that
$$\omega(\xi\G)=\frac{\mu_H((H\cap {g_\xi}^{-1}\G g_\xi)\ba H)}{\mu_G(\G\ba G)}>c\quad\text{ for all $\xi\in \cup_m\Cal O_m$.}$$ Since $\omega (\Cal
O_m) \ge h_m \cdot c$ where $h_m$ is the number of disjoint $\G$-orbits
in $\Cal O_m$, we may also assume that $h_m$ is constant,
say $r$, for all $m$, by passing
to a subsequence. Now write $\Cal O_m=\cup_{i=1}^r\xi_{m_i}\G$.
It suffices to show that for each $1\le i\le r$,
$\pr(\xi_{m_i}\G)$ is the same set for infinitely many $m$. For,
this implies that $\pr(\Cal O_m)$ is the same set
for infinitely many $m$, which contradicts (1).
 Fix $1\le i\le r $ and set $\xi_{m_i}=\xi_m$ for simplicity.

 It follows from Propositions  \ref{supp0} and \ref{Non-div} that
 there exists a compact
subset $C$ of $\G\ba G$ such that
$\G\ba \G g_{\xi_m} H\cap C \ne\emptyset$
for all $m\in \n$.
Hence we may choose $g_{\xi_m}$ so that $\{g_{\xi_m}: m\in \n\}$ is
relatively compact.

On the other hand, if $\delta_g$ denotes the factor by which the
volumes of subsets gets multiplied under the transformation $h \to
ghg^{-1}$, $h\in H$, then
$$\mathcal V ((H_{\xi_m}\cap \G )\ba H_{\xi_m})=
\delta_{g_{\xi_m}}\cdot \omega(\xi_m \G) $$
up a uniform constant multiple depending only
on the choice of Haar measure
 $\mu_H$.
Since $\{g_{\xi_m}: m\in \n\}$ is relatively compact,
$\sup_m \delta_{g_{\xi_m}} <\infty$.

Therefore
 $$\sup_m\mathcal V
((H_{\xi_m}\cap \G )\ba H_{\xi_m}) \le \sup_m \delta_{g_{\xi_m}}\cdot \sup
_m\omega(\Cal O_m) <\infty .$$
By Theorem \ref{fv}, this implies that $H_{\xi_m}\cap \G$ are all equal
to each other by passing to a subsequence.
Since $H_{\xi_m}\cap \G$ is Zariski dense in $H_{\xi_m}$ by Borel
density theorem, it follows that $H_{\xi_{m}}$ are all
equal to each other, that is, $g_{\xi_m}^{-1}g_{\xi_k}\in N(H)$
for all $m,k$. Since $[N(H):H]<\infty$, by passing to a subsequence,
 we have $$g_{\xi_m}^{-1}g_{\xi_k}\in H, \quad\text{ and hence} \quad
\pr(\xi_m\G)=\pr(\xi_k \G)$$
for all $m,k$.
 This finishes the proof.
\end{proof}

Observe that
for compact subsets $\Omega\subset V_1 $, $C_0\subset G$ such that
$H\ba HC_0=\Omega$ and $\xi\in V_m(\z)$,
\begin{equation}\label{ifo} N(\xi\G,  \Omega) = 0\quad \text{if
and only if } \quad \G\ba \G g_{\xi} H\cap \G\ba \G C_0^{-1} = \emptyset .
\end{equation}

\begin{Prop}\label{proper} Consider a sequence $\{\xi_m\in V_m(\z)\}$
such that $\{\pr(\xi_m)\}$ is relatively compact in $V_1$.
 Suppose that $L$ is a closed
subgroup of $G$ containing
$ H_{\xi_m}^0$ for all $m$.
such that $L\cap \G$ is a lattice in $L$.
 Then for any compact subset $\Omega$ of
$ V_1$, there exists a finite subset $\Lambda_{\Omega}\subset \G$
 such that for all $m$,
$$ \operatorname{pr} (\xi_m) \G
\cap \Omega \subset \operatorname{pr} (\xi_m) (L\cap \G)
 \Lambda_{\Omega}.$$
\end{Prop}
\begin{proof} Let $\Omega_0$ be a compact subset of $G$ such that
$H\ba H\Omega_0=\Omega$.
Write $H$ as a disjoint union $\cup_{i=1}^k h_iH^0$.
By the assumption, there exists a choice of
$\{g_{\xi_m}\}$ so that $\{g_{\xi_m}\in G\}$ is
relatively compact.
Let $\Omega_1\subset G$ be a compact
subset which contains $\{g_{\xi_m}h_i: 1\le i\le k,\,
m\in \n\}\Omega_0$.

Since $L\cap \G$ is a lattice in $L$, $\G\ba \G L$ is closed
in $\G\ba G$ [Rag], and this implies easily that $L\ba L\G$ is
closed in $L\ba G$. Since $L\ba L\G$ is a closed countable
subset of $L\ba G$,  it follows from Baire category theorem
that there exists at least one isolated point.
Since $\G$ acts transitively on $L\ba L\G$,
every point of $L\ba L\G$ is an isolated point.
 Therefore $L\ba
L\G$ is discrete in $L\ba G$.
Hence there exists a finite
subset $\Lambda_{\Omega}$ of $ \G$
such that
$$L\Omega_1 \cap L\G\subset L \Lambda_{\Omega} .$$
Note that $$g_{\xi_m}H\Omega_0 \cap g_{\xi_m}H^0g_{\xi_m}^{-1}\G
\subset Lg_{\xi_m} (\cup_{i=1}^k h_i)
\Omega_0 \cap L\G\subset L\Omega_1 \cap
L\G\subset L \Lambda_{\Omega}.$$

Hence
\begin{equation}\label{ho} H\Omega_0 \cap H^0
g_{\xi_m}^{-1}\G \subset g_{\xi_m}^{-1}L \Lambda_{\Omega}.
\end{equation}

If $v_0 x\in \pr (\xi_m \G)\cap \Omega$ for $x\in G$,
then
$$x=hg_{\xi_m}^{-1}\gamma=h'w
$$ for some $h, h'\in H$,
$\gamma\in \G$ and $w\in \Omega_0$.

Then $$g_{\xi_m}^{-1}\gamma \in H\Omega_0 \cap
H^0g_{\xi_m}^{-1}\G,$$ and hence by (\ref{ho}),
$$g_{\xi_m}^{-1}\gamma=g_{\xi_m}^{-1}g\gamma_1 \quad\text{
for some
$g\in L$ and $\gamma_1\in \Lambda_{\Omega}$.}$$ In
particular, $g=\gamma \gamma_1^{-1}\in L\cap \G$.
 Therefore
$$x=hg_{\xi_m}^{-1}g\gamma _1\in
 H g_{\xi_m}^{-1}(L\cap \G)\Lambda_{\Omega}$$
proving
 $$  \operatorname{pr} (\xi_m \G)\cap \Omega \subset \operatorname{pr} (\xi_m) (L\cap \G)
 \Lambda_{\Omega}.$$
\end{proof}

\section{Asymptotic behavior of $\Cal O_m$}\label{last}
We start by recalling the following
theorem of Dani and Margulis.
 For any two closed subgroups $U$ and $L$ of a connected Lie group $G$,
set $$X(U, L):=\{g\in G: gU\subset Lg\}.$$
\begin {Thm} [DM2, Theorem 3]\label{dm} Let $G$ be a connected Lie group and $\G$ a lattice in $G$.
Let $U=\{u(t)\}$ be a unipotent one-parameter subgroup of $G$ and
let $\psi$ be a bounded continuous function of $\gg$. Let $K$ be a compact
subset of $\gg$ and let $\e >0$ be given. Then
there exist finitely many proper closed subgroups $L_1, \cdots , L_k$
such that
$L_i\cap \G$ is a lattice in $L_i$ for each $1\le i\le k$,
 and compact subsets $C_1,
\cdots , C_k$ of $X(U, L_1), \cdots , X(U, L_k)$
respectively, for which the following holds:
for any compact subset $F \subset K -\cup_{i=1}^k\G \ba \G C_i$, there exists
 $T_0 \ge 0$ such that for all $x\in F$ and $T>T_0$,
\begin{equation*}
\left| \frac{1}{T}\int_{0}^T\psi(xu(t))\, dt -\int_{\gg}\psi \,
d\mu_G \right|\le \e .\end{equation*}
\end{Thm}

In fact, it is shown in the proof of the above theorem [DM2] that
the subgroups $L_i$ can be taken so that $\operatorname{Ad} (L_i\cap \G)$ is Zariski
dense in $\operatorname{Ad} (L_i)$
 as well where $\operatorname{Ad} $ denotes the
 adjoint representation of $G$.

We keep the same notation from section 2 for
 $G$, $\rho:G \to \GL(V)$, $\G$, $v_0$, etc.
Let $\Cal O_m$ be a $\G$-invariant subset of $V_m(\z)$.
As before, we assume that the subgroup $H$, which is the stabilizer
of $v_0$, is a semisimple real algebraic
subgroup of $G$ without compact factors.
   Set $X=\G\ba G$. We assume without loss of generality
$\mu_G(X)=1$. For each $\xi\in V_m(\z)$, we denote by $\nu_{\xi\G}$
the unique $H$-invariant measure on $X$ supported on
$\G\ba \G g_\xi  H$ with the total measure given by $\omega (\xi\G)$.

Denote by $\pi$ the canonical projection from $H^0\ba G$ to $H\ba G$.
By the identification of $V_1=v_0G$ with $H\ba G$,
we consider $\pr(\Cal O_m)$ as a subset of $H\ba G$.
Set $\frak O_m=\pi^{-1}(\pr(\Cal O_m))$.
Note that ${\frak O_m}$ is $\G$-invariant and has finitely
many $\G$-orbits. For each $\eta \G\subset \frak O_m$,
the notation $\nu_{\eta\G}$ denotes the $H^0$-invariant measure
 on $\G\ba G$ supported on $\G\ba \G g_{\eta} H^0$ with the total measure
given by $$\omega(\eta \G):=\frac{\mu_{H^0} ((H^0\cap
g_{\eta}^{-1}\G g_{\eta})\ba H^0)}{\mu_G(\G\ba G)}$$ where $g_\eta$
is any element in $G$ such that $H^0\ba H^0g_\eta=\eta$. Here
$\mu_{H^0}$ is simply the restriction of
 $\mu_H$ to $H^0$.
Note that if $\pi(\eta_1\G)=\pi(\eta_2\G)$,
then $\omega(\eta_1\G)=\omega(\eta_2\G)$ and $\omega(\eta\G)\le
\omega(\xi\G)$ if $\pi(\eta)=\pr(\xi)$.

For $\xi\G\subset \Cal O_m$, it is not hard to check that
$$\nu_{\xi\G}=\sum_{\eta \G\subset \pi^{-1}(\pr(\xi\G))}
 \nu_{\eta\G}$$
where the sum is taken over the disjoint $\G$-orbits $\eta \G$
in $\frak O_m$
such that
$\pi(\eta\G)=\pr(\xi\G)$. The number such $\G$-orbits is clearly bounded by
$[H:H^0]$.

Therefore
\begin{equation}\label{oo}\sum_{\xi\G\subset \Cal O_m}\nu_{\xi\G}=
 \sum_{\eta \G\subset \frak O_m}
 \nu_{\eta \G} \text{ and }\omega(\Cal O_m)=\sum_{\eta \G\subset \frak O_m}
 \omega(\eta \G).
 \end{equation}

 Let $U=\{u(t)\}$ be a unipotent one-parameter
subgroup of $H^0$ not contained in any proper closed normal subgroup
of $H^0$.
For given compact
subset $K\subset X$, $\e >0$ and a bounded continuous function
$\psi$ on $X$, let $L_i$ and $C_i$, $1\le i\le k$ be as in the
above theorem, with respect to the given triples $(K, \e, \psi)$.
Set $$\Cal S(K, \psi, \e):=K\cap (\cup_{i=1}^k \G\ba \G C_i)$$
and let $\Cal G (K, \psi, \e)$ denote the complement of
$\Cal S(K, \psi, \e)$ inside $K$.

\begin{Prop}\label{di2} Fix a compact subset $K$ of $ X$ with a non-empty
interior.
Suppose that for any $\e >0$ and for any continuous
 function $\psi$ on $X$ with
compact support,
$$ \lim_{m \to \infty}
\frac{\sum_{\eta\G\subset \frak O_m} \{\omega (\eta\G): \G\ba
\G g_\eta H^0\cap \Cal G (K,
\psi, \e)\ne \emptyset \}}
{\omega (\Cal O_m)}=1 . $$
Then
$$ \frac{1}{\omega(\Cal O_m)}\sum_{\xi\G\subset\Cal O_m}
\nu_{\xi\G}
\to \mu_{ G}\quad \text{  as $m \to \infty$}.$$
\end{Prop}
\begin{proof}
We set
$$A_m(\e):=\{\G\ba \G g_\eta H^0 :\eta\G\subset \frak O_m,
\G \ba \G g_\eta H^0\cap K \subset \Cal S (K, \psi ,
 \e)  \}$$
and $$B_m( \e ):= \{\G\ba \G g_\eta H^0
: \G\ba \G g_\eta H^0 \cap \Cal G(K, \psi,
 \e)\ne \emptyset\}.$$

The assumption implies that
\begin{equation}\label{am}
 \lim_{m \to \infty} \frac{ \sum \{\omega (\eta\G)
 : \G\ba \G g_\eta H^0 \in
A_m(\e)\} }{\omega(\Cal O_m) }=0 ;\quad
  \lim_{m \to \infty} \frac{  \sum \{\omega (\eta\G) : \G\ba \G g_\eta H^0 \in
B_m(\e)\} }{\omega(\Cal O_m)}=1 .\end{equation}

 By the Birkhoff ergodic theorem and Moore's ergodicity theorem,
$U$ acts ergodically
with respect to each $\frac{1}{\omega(\eta\G)}
\nu_{\eta\G}$ and the following subset $R$
has the zero co-measure in $H^0$:
$$R=\{ h\in H^0: \G\ba \Gamma g_\eta hu(t)
\text{ is uniformly distributed in $\G\ba \Gamma g_\eta H^0$ w. r.
t. $\frac{1}{\omega(\eta\G)}\nu_{\eta\G}$} \}.$$

Let $\G\ba \G g_\eta H^0 \in B_m(\e)$.
Since
 $\Cal G (K, \psi, \e)$ is open in $K$ and
 $R$  has co-measure $0$ in $H^0$,
we may assume by a suitable choice for
$g_\eta$ that $\G\ba \G g_\eta \in \Cal G (K, \psi, \e)$ and
$$\lim_{T\to \infty}\frac{1}{T}
\int_0^T\psi(\G\ba \G g_\eta u(t))dt=
 \frac{1}{\omega(\eta\G)} \int_{X} \psi \, d \nu _{\eta\G}.$$
Therefore by
applying Theorem \ref{dm} to each singleton $F=\{\G\ba \G g_\eta\}$,
we obtain that
for any $\G\ba \G g_\eta H^0 \in B_m(\e)$,
\begin{equation}\label{dme}
\left|\frac{1}{\omega(\eta\G)} \int_{X} \psi \, d \nu _{\eta\G}
-\int_{X}\psi \, d\mu_G\right| \le \e .
\end{equation}

Now
\begin{multline*}
\left| \sum_{\eta\G\subset\frak O_m} \int_{X}\psi \, d\nu_{\eta\G}
-\int_X \psi \,d\mu_G \right| \le
 \sum_{\G\ba \G g_\eta H^0\in A_m(\e)} \left|\int_{X}\psi \,
d\nu_{\eta\G}  -\omega(\eta\G) \int_X \psi \, d \mu_G\right|  \\ +
\sum_{\G\ba \G g_\eta H^0\in B_m(\e)}
 \left|\int_{X}\psi \, d\nu_{\eta\G} -\omega(\eta\G)
\int_X \psi\, d \mu_G \right|
\end{multline*}
By (\ref{dme}), the above is again less than or equal to
$$ (\|\psi\|_{\infty}+\|\psi\|_1)
 ( {\sum\{ \omega (\eta\G) :\G\ba
\G g_\eta H^0\in A_m(\e)\} }) +{ \e}( \sum\{ \omega (\eta \G) :\G\ba
\G g_\eta H^0\in B_m(\e)\} ).$$

By applying Lebesgue dominated convergence theorem,
we deduce from (\ref{am})
$$\limsup_{m \to \infty} \left|\frac{1}{\omega (\Cal O_m)}
\sum_{\eta\G\subset\frak O_m}
\left(\int_{X}\psi \, d\nu_{\eta\G}\right)
-\int \psi d \mu_G \right| \le \e .$$
Since $\e >0$ is arbitrary, we have
for any bounded continuous function
$\psi$ on $X$,
 $$\lim_{m \to \infty}\frac{1}{\omega(\Cal O_m)}
\sum_{\eta\G\subset\frak O_m}
\left(\int_{X}\psi \, d\nu_{\eta\G}\right)
=\int \psi \,d \mu_G .$$

This proves our claim by (\ref{oo}).
\end{proof}

\begin{Lem}[EMS2, Lemma 5.1] \label{EMS} Let $G$, $H$ and $L$ be
 connected real algebraic groups such that $H\subset L\subset G$.
 If at least one of $G$, $H$, and $L$ is reductive, then $X(H,
L)$ is a union of finitely many closed double cosets of the form
$L \cdot g\cdot Z(H)$ where $g\in X(H, L)$.
\end{Lem}

 \noindent {\bf Proof of Theorem \ref{Main}} Since $G^0$ has
a finite index in $G$ and $V_1$ consists of finitely many open $G^0$
orbits, it suffices to prove the theorem for each $G^0$-orbit. Hence
we may assume that $G$ is connected without loss of generality.
Since
 $\{\Cal O_m\}$ is not focused, there exists an arithmetic subgroup
$\G\subset G(\q)$ which preserves each $\Cal O_m$ and $\{\Cal O_m
\}$ is not focused with respect to $\G$.

First, for some compact subset $C$ of $V_1$,
 $$\limsup_{m \to \infty}\frac{\sum_{\xi\G\subset\Cal O_m}
\{\omega(\xi\G): \operatorname{pr} (\xi)\G\cap C =\emptyset
\}}{\omega(\Cal O_m)} =0 . $$ Note that the same holds for any compact
subset of $V_1$ containing $C$.

Hence it follows from the observation (\ref{ifo}) that for some relatively compact open subset $C_0$
of $G$, we have
\begin{equation}\label{lims}\limsup_{m \to \infty}
\frac{\sum_{\eta\G\subset\frak O_m}
\{\omega(\eta\G): \G\ba
 \G g_\eta H^0 \cap \G\ba \G C_0 =\emptyset \}}{\omega (\Cal O_m)} =0 .
 \end{equation}

Set $K=\G\ba (\G \overline {C_0})$ and $K'=\G\ba (\G {C_0})$. By
Propositions \ref{Av},
 and \ref{di2}, it suffices to show that
 for any $\e >0$ and for any bounded continuous
 function $\psi$ on $X$,
\begin{equation*} \lim_{m \to \infty}
\frac{\sum_{\eta\G\subset\frak O_m} \{\omega(\eta\G): \G\ba \G
g_\eta H^0\cap \Cal G (K, \psi, \e)\ne \emptyset \}} {\omega(\Cal
O_m)}=1 . \end{equation*} Suppose not. Since the orbits $\G\ba \G
g_\eta H^0$ disjoint from $K'$ can be ignored by (\ref{lims}),
 it follows that there exist
a bounded continuous
function $\psi$ on $X$ and an $\e >0$ such that
 \begin{equation}\label{po}\limsup_{m \to \infty}
\frac{\sum_{\eta\G\subset \frak O_m} \{\omega(\eta\G):
\emptyset\ne K'\cap \G\ba \G g_\eta H^0 \subset \Cal S (K,
\psi, \e) \} }{\omega(\Cal O_m)}> 0\end{equation}
Let $L_i$ and $C_i$, $1\le i\le k$ be the
subgroups and compact subsets in $X(U,L_i)$ respectively,
 used in the definition of
 $\Cal S (K, \psi, \e)$.
Since $\Cal S(K, \psi, \e)$
is contained in the
finite union $\cup_{i=1}^k \G\ba \G C_i$,
 there exists $1\le i\le k$
such that (\ref{po}) holds with $\G\ba \G C_i$ in place of
 $\Cal S(K, \psi, \e)$. Without loss of generality, we assume $i=1$.
By Lemma \ref{EMS}, there exists $g\in G$ such that
$gUg^{-1}\subset L_1$ and
\begin{equation}\label{po1}\limsup_{m \to \infty}
\frac{\sum_{\eta\G\subset \frak O_m} \{\omega(\eta\G): \emptyset \ne
\G g_\eta H^0 \cap K '\subset  \G (C_1 \cap L_1 gZ(U)) \}
 }{\omega(\Cal O_m)}> 0. \end{equation}

Whenever $$
\emptyset \ne \G g_\eta H^0 \cap K '\subset  \G (C_1\cap L_1 gZ(U))$$
 we may assume
that $g_\eta \in C_1\cap L_1gZ(U)$ by replacing $\eta$ and $g_\eta$
by suitable elements  in $\pi^{-1}(\pr(\xi\G))$ and $g_\eta H^0$ respectively.
We may also assume
that the set $\{h\in H^0: g_\eta h\in C_1\cap L_1 gZ(U)\}$ has
a positive measure, since $K'$ is open. Note that
$g_\eta h\in C_1\cap L_1 gZ(U)$ implies that $h Uh^{-1}
\subset {g_\eta}^{-1} L_1 g_\eta$. By a similar argument as in
the proof of Proposition \ref{Non-div},
it follows that $H_{\xi}^0\subset L_1$, i.e., $g_\xi \in X(H^0,L_1)$,
whenever
$\emptyset \ne \G g_\eta H^0 \cap K '\subset  \G (C_1\cap L_1 gZ(U))$
and $\pi(\eta)\in \pr(\xi\G)$.
Applying Lemma \ref{EMS} again,
 we deduce from (\ref{po1}) that for some
 $g_0\in G$ such that $g_0H^0g_0^{-1}\subset L_1$,
 \begin{equation*}\limsup_{m \to \infty}
\frac{\sum_{\eta\G\subset \frak O_m} \{\omega(\eta\G):
 g_\eta \in C_1\cap L_1g_0Z(H^0)\} }{\omega(\Cal O_m)}> 0 .\end{equation*}

 Since $g_\eta\in C_1\cap L_1 g_0Z(H)$ and $\pi(\eta)=\pr(\xi)$
 implies that $g_\eta H^0g_\eta^{-1}=H_{\xi}^0\subset
L_1$,
it follows from Proposition \ref{proper} that
 for any compact subset
$\Omega$ of $ V_1$, there exists
a finite subset $\Lambda_{\Omega}\subset \G$
such that for all $g_\eta\in C_1\cap L_1 g_0Z(H)$,
$$\pi(\eta) \G\cap \Omega\subset
(H {g_\eta}^{-1}) (L_1\cap \G)
\Lambda_{\Omega} \subset
H{( L_1g_0 Z(H)\cap C_0)}^{-1} (L_1\cap \G)
\Lambda_{\Omega} .$$
Hence we have shown that for any compact subset $\Omega\subset V_1$,
 \begin{equation*}\limsup_{m \to \infty}
\frac{\sum_{\eta\G\subset\frak O_m}
 \{\omega(\eta\G): \pi(\eta) \G\cap \Omega\subset
 H{( Z(H)g_0^{-1} L_1\cap C_0 ^{-1})} (L_1\cap \G)
\Lambda_{\Omega}\} }{\omega (\Cal O_m)}> 0 ,\end{equation*}
for some finite subset $\Lambda_{\Omega}\subset \G$.
Hence for some $\gamma_\Omega\in \G$,
\begin{equation*}\limsup_{m \to \infty}
\frac{\sum_{\xi\G\subset\Cal O_m}
 \{\omega(\xi\G): \pr(\xi) \G\cap \Omega\subset
 v_0{( Z(H)g_0^{-1} L_1\cap C_0 ^{-1})} (L_1\cap \G)
\gamma_{\Omega}\} }{\omega (\Cal O_m)}> 0 ,\end{equation*}
By the remark following
Theorem \ref{dm},
$L_1\cap \G$ is a Zariski dense lattice in $L_i$.

Hence the sequence $\{\Cal O_m\}$ is
 focused, yielding
contradiction. This finishes the proof.

\noindent{\bf Proof of Theorem \ref{maximal}} Without loss of
generality we may assume that $G$ is connected. Fix any $\G\in \Cal
A_V$ preserving each $\Cal O_m$. If (\ref{p1}) does not hold, then $\{\Cal
O_m\}$ is focused by Theorem \ref{main}. Since $H^0$ is a maximal
connected closed subgroup of $G$, it follows that there exists $g\in
G$ such that for any compact subset $\Omega \subset V_1$
 \begin{equation*}\limsup_{m \to \infty}
\frac{\sum_{\xi\G\subset\Cal O_m} \{\omega(\xi\G): \pr(\xi)\G\cap
\Omega \subset v_0Z(H^0)H^0g \G
 \} }{\omega(\Cal O_m)}> 0 .\end{equation*}

We claim that
\begin{equation}\label{mf}
\limsup_{m \to \infty} \frac{\sum_{\xi\G\subset\Cal O_m}
\{\omega(\xi\G): \pr(\xi)\G\cap  v_0Z(H^0)H^0g \G
\ne \emptyset \}}{\omega(\Cal O_m)} > 0 \end{equation}
Suppose not. Then
it follows that
 for any compact subset $\Omega$ of $V_1$,
 $$\limsup_{m \to \infty}\frac{\sum_{\xi\G\subset\Cal O_m}
\{\omega(\xi\G): \operatorname{pr} (\xi)\G\cap \Omega =\emptyset
\}}{\omega(\Cal O_m)} >0 . $$
This is equivalent to saying that
 for any compact subset $C$ of $\G\ba G$,
 $$\limsup_{m \to \infty}\frac{\sum_{\eta\G\subset\frak O_m}
\{\omega(\eta\G): \G\ba \G g_\eta H^0\cap C =\emptyset
\}}{\omega(\Cal O_m)} >0 . $$
Hence there exists a sequence $\eta_{m}\in \frak O_m$
such that for any compact subset $C\subset X$,
there exists $m$ such that $ \G\ba \G g_{\eta_m} H^0\cap C=\emptyset$.
On the other hand, since $H^0$ is a proper maximal
closed subgroup of $G^0$, by Proposition \ref{Non-div},
any weak limit of $\{\nu_{\eta_m\G}\}$ in $\Cal P(X\cup\{\infty\})$
 is supported on $X$. This is a contradiction by Proposition \ref{supp0}.
Hence (\ref{mf}) is proved.

It is easy to check that
 $\pr(\xi)\G\cap v_0Z(H^0)H^0 g\G \ne\emptyset$ implies
$ \gamma H_{\xi}\gamma ^{-1}=gHg^{-1}$
for some $\gamma\in \G$.
Hence (\ref{mf}) implies:
\begin{equation*}
\limsup_{m \to \infty} \frac{\sum_{\xi\G\subset\Cal O_m}
\{\omega(\xi\G):  \gamma H_{\xi}\gamma ^{-1}=gHg^{-1}
\text{ for some $\gamma\in \G$}
 \}}{\omega(\Cal O_m)} > 0 \end{equation*}

It follows that
 there exists $\xi_0\in \Cal O_{m_0}$ for some $m_0$ such that
\begin{equation*}
\limsup_{m \to \infty} \frac{\sum_{\xi\G\subset\Cal O_m}
\{\omega(\xi\G):  \gamma H_{\xi}\gamma ^{-1}=H_{\xi_0}
\text{ for some $\gamma\in \G$}
 \}}{\omega(\Cal O_m)} > 0 \end{equation*}

Observe that the condition
$ \gamma H_{\xi}\gamma ^{-1}=H_{\xi_0}$ implies that
$g_{\xi_0}^{-1} \gamma g_{\xi}\in N(H)$,
and the condition
$g_{\xi_0}^{-1} \gamma g_{\xi}\in H$ implies
that $\pr(\xi)\G=\pr(\xi_0)\G$.

Since $H$ has a finite index in the normalizer
$N(H)$, it follows that
\begin{equation}\label{mmm}
\limsup_{m \to \infty}\frac{\sum_{\xi\G\subset\Cal O_m}
 \{\omega(\xi\G) :
\pr (\xi)\G =\pr (\xi_{0})\G)\}}{
\omega (\Cal O_m)} >0 .\end{equation}

Note that if $\pr (\xi)\G=\pr (\xi_{0})\G$, then
$\omega(\xi\G)=\omega(\xi_0\G)$.
Since there can be at most one $\G$-orbit $\xi\G$ in $\Cal O_m$ such that
$\pr (\xi)\G=\pr (\xi_{0})\G$,
 (\ref{mmm}) implies that
\begin{equation*}
\omega (\xi_0\G) \cdot
\limsup_{m \to \infty}\frac{1}{
\omega(\Cal O_m)} >0 \quad \text{ or
equivalently }\quad
\liminf_{m \to \infty}\omega(\Cal O_m) <\infty .\end{equation*}
By Lemma \ref{q}, this
 contradicts the assumption on $\{\Cal O_m\}$.
 Hence the proof is now complete.

\section{Examples}\label{example}
\begin{Thm}\label{quad}
Let $Q$ be an integral quadratic form of signature $(r, s)$
where $r+s\geq 4$, $r\ge 2$ and $ s\ge 1$.
 For any compact subset $\Omega$ of
$ V_1$ with boundary of measure $0$,
 $$N_m(Q, \Omega)
\sim
 \omega_m \cdot \operatorname{vol}{(\Omega)}$$
as $m \to \infty$ along the fundamental discriminants.
 \end{Thm}
\begin{proof} Consider the
standard representation of the orthogonal
group $\operatorname{O}(Q)$ on
$V:=\br^{r+s}$. Let $V_m$ be the level set $\{x\in V: Q(x)=m\}$.
 By Witt's theorem the orthogonal group
$\operatorname{O}(Q)$ acts transitively on each $V_1$.
The stabilizer of a vector $v_0$ in $V_{1}$
 is isomorphic to $\operatorname{O}(r-1, s)$.
 Note that the assumptions on the size
of the parameters $r$ and $s$ guarantee that
$\operatorname{O}(r-1,s)$ is non-compact
and simple. It is well known that $\operatorname{O}(r-1,s)^0$ is a
maximal connected subgroup of $\operatorname{O}(r,s)^0$.
If we set $\G:=\operatorname{O}(Q)\cap \operatorname{SL}_{r+s}(\z)$,
then $\G\in\Cal A_V$.
Under the assumption that $r+s\ge 4$ and $rs=1$,
$V_m(\z)\ne \emptyset$ for all fundamental discriminants $m$ (see
[Oh]). Hence we may take $\Cal O_m=V_m(\z)$ to apply theorem \ref{maximal}.
To check the condition (\ref{prim}), note that if $m\ne k$ are
fundamental discriminants, then $\sqrt m ^{-1} V_m(\z)\cap \sqrt k
^{-1} V_k(\z)= \emptyset$; otherwise this would imply that $\sqrt{
(m/k)} \in \q$, which can be seen to be false by an easy
computation. Therefore
 Theorem \ref{maximal} implies the claim.
\end{proof}

Let $V:=\{X\in \operatorname{M}_{2n}(\br): X^t=-X\}$ be the space
 of skew-symmetric matrices,
so that the Pfaffian on $V$ is defined by
$$\operatorname{Pff}^2(X)= {\operatorname{Det}(X)}\quad\text{and}
 \quad \operatorname{Pff}(v_0)=1$$ where
$$v_0=\begin{pmatrix} 0& I_n\\-I_n&  0\end{pmatrix} .$$
E.g., for $n=2$, we have $\text{Pff}(x_1, \cdots , x_6)
=x_1x_2-x_3x_4+x_5x_6$.

\begin{Thm}\label{pff} Let $n \ge 2$.
For any compact subset $\Omega$ of $V_1$ with boundary of measure $0$,
\begin{equation*}
N_m(\operatorname{Pff}, \Omega) \sim_{m\to \infty} \omega_m \cdot
\operatorname{vol}{(\Omega)}.
\end{equation*}
\end{Thm}
 \begin{proof} Consider the representation
$\rho:\SL_{2n}(\br) \to \GL(V)$ defined by
$$\rho(A)(X)=A^{t}XA$$
where $A\in \SL_{2n}(\br)$ and $X\in V$.
 It is well known that $\SL_{2n}(\br)$
acts transitively on $V_1=\{X\in V: \operatorname{Pff}(X)=1\}$, so that
 $V_m=m^{1/n} v_0 \SL_{2n}(\br)$.
The stabilizer $H$ of $v_0$ is the symplectic group
$\Sp_{2n}(\br)$ corresponding to $v_0$, which is a maximal connected
closed subgroup of $\SL_{2n}(\br)$.
Clearly $V_m(\z)\ne\emptyset$ and $\SL_{2n}(\z)$ preserves
$V_m(\z)$ for each $m$.

To check the condition \ref{prim} of Theorem \ref{maximal},
 suppose for a given $m_0\in \n$
 that $\pr (V_{m_0}(\z))=\pr (V_{m}(\z))$.
Then  $$m^{-1/n}\begin{pmatrix} & & m &0 \\ && 0 & I_{n-1}   \\
 -m & 0 && \\ 0 &-I_{n-1} &&
\end{pmatrix}
\in m_0^{-1/n} V(\z).$$ Hence $m^{-1/n}=m_0^{-1/n} k$ for some
$k\in \n$. This leads to $m=m_0k^{-n} \le m_0$. Hence
$$|\{m\in \n: \pr (V_{m_0}(\z))=\pr (V_{m}(\z))\}|\le m_0 <\infty .$$
 Hence the claim follows from Theorem \ref{maximal}.
 \end{proof}

Since the explicit representatives of $\SL_{2n}(\z)$-orbits on
$V_m(\z)$ can easily be written down in the Pff case, we can
compute $\omega_m$, using the local density formula given in [GY].
For instance, for a prime $p$,
$$\omega_p=C \cdot \sum_{i=0}^{2n-2}p^i$$
 for some constant $C$ independent of
$p$.

 \begin{Thm}\label{det} Let $n \ge 2$.
For any compact subset $\Omega$ of $V_1=\SL_n(\br)$  with boundary of measure $0$,
\begin{equation*}
N_m(\operatorname{Det}, \Omega) \sim_{m\to \infty} \omega_m \cdot
\operatorname{vol}{(\Omega)}.
\end{equation*}
\end{Thm}

\begin{proof} Consider the representation $\SL_n(\br) \times \SL_n(\br)$ on
the space $V=M_n(\br)$ given by $X(A,B)=AXB^{-1}$ where $X\in V$
and $A, B\in \SL_n(\br)$.
The stabilizer of $I_n$ is given by the diagonal embedding of $\SL_n(\br)$
in the product $\SL_n(\br)\times \SL_n(\br)$, which
is a maximal connected closed subgroup. Clearly $V_m(\z)\ne \emptyset$ and
$\SL_n(\z)\times \SL_n(\z)$ preserves $V(\z)$.
 The condition \ref{prim} can be checked
similarly as in the case of Pfaffian.
Hence our claim Theorem \ref{maximal}.
\end{proof}

For the determinant case, the constant $\omega_m$ is well known
from the theory of Hecke operators. (cf. [Sa], or [GO]). For
instance, for any fixed $k\in \n$ and a prime $p$,
 $$\omega_{p^k}\sim_{p\to \infty} C \cdot p^{k(n-1)}.$$

\enddocument